\documentclass[journal,12pt,onecolumn,draftclsnofoot,]{IEEEtran}
\usepackage[a4paper, total={7in, 10in}]{geometry}
\usepackage{blindtext}
\usepackage{graphicx}
\usepackage{multirow}
\usepackage{amsmath,amsthm}
\usepackage{array} 
\usepackage{etoolbox}
\usepackage{multicol}
\usepackage{multirow}
\usepackage{hhline}
\usepackage{ifthen}
\usepackage{array}
\usepackage[noprefix]{nomencl}
\usepackage{nomencl}
\usepackage{siunitx}
\usepackage{algorithm} 
\usepackage{algpseudocode}
\usepackage{pifont}
\usepackage{subcaption} 
\usepackage{url}
\makenomenclature
\renewcommand{\nomgroup}[1]{
\ifthenelse{\equal{#1}{A}}{\vspace{0.3cm} \item[\bfseries Indexes and Sets]}{
\ifthenelse{\equal{#1}{B}}{\vspace{0.3cm} \item[\bfseries Sets]}{
\ifthenelse{\equal{#1}{C}}{\vspace{0.3cm} \item[\bfseries Parameters]}{
\ifthenelse{\equal{#1}{D}}{\vspace{0.3cm} \item[\bfseries First-stage Variables]}{
\ifthenelse{\equal{#1}{E}}{\vspace{0.3cm} \item[\bfseries Second-stage Variables]}{}
}}}}}

\usepackage{ifthen}
\usepackage{float}
\usepackage{hhline}
\usepackage{tikz}
\usepackage{pgfplots}
\pgfplotsset{width=9cm,compat=1.3}
\usepackage{xcolor}
\usepackage{pgfplotstable}
\usepackage{cite}
\usetikzlibrary{patterns}

\definecolor{ref}{rgb}{0.65,0.65,0.65} %{0.4,0.8,0.85}
\definecolor{lhmm}{rgb}{0.9,0.6,0.5}
\definecolor{ghmm}{rgb}{0.7,0.9,0.35}
\definecolor{lsvm}{rgb}{0.9,0.8,0.25}
\definecolor{gsvm}{rgb}{0.4,0.8,0.9}

\usetikzlibrary{arrows,automata,shapes,chains,shadows,calc}
\usepackage{caption}
\usepackage{subcaption}
\usepackage{comment}
\usepackage{hhline,multirow,float} 
\usepackage{algorithm}
\usepackage{epstopdf}
\ifCLASSINFOpdf
\else
\fi
\begin{document}
\title{An Efficient Robust Solution to the Two-Stage Stochastic Unit Commitment Problem}
\author{Ignacio~Blanco~and~Juan~M.~Morales
\thanks{Ignacio Blanco and J. M. Morales (corresponding author) are with the Technical University of Denmark, DK-2800 Kgs. Lyngby, Denmark (email addresses: \{igbl, jmmgo\}@dtu.dk), and their work is partly funded by DSF (Det Strategiske Forskningsr{\aa}d) through the CITIES research center (no. 1035-00027B).}}
\maketitle

\begin{abstract}
This paper provides a reformulation of the scenario-based two-stage unit commitment problem under uncertainty that allows finding unit-commitment plans that perform reasonably well both in expectation and for the worst case. The proposed reformulation is based on partitioning the sample space of the uncertain factors by clustering the scenarios that approximate their probability distributions. The degree of conservatism of the resulting unit-commitment plan (that is, how close it is to the one provided by a purely robust or stochastic unit-commitment formulation) is controlled by the number of partitions into which the said sample space is split. To efficiently solve the proposed reformulation of the unit-commitment problem under uncertainty, we develop a parallelization and decomposition scheme that runs as a column-and-constraint generation procedure. Finally, we analyze the quality of the solutions provided by this reformulation for a case study based on the IEEE 14-node power system and test the effectiveness of the proposed parallelization and decomposition solution approach on the larger IEEE 3-Area RTS-96 power system.   

\end{abstract}

% \fi
\IEEEpeerreviewmaketitle

%The main notation used throughout the paper is stated below for quick reference. Other symbols are defined as required.
\section*{Nomenclature}
%Notation
The notation used throughout the paper is stated below for quick reference. Other symbols are defined as required.
\subsection{Indexes and Sets}
\begin{IEEEdescription}
\item[$T$] Set of time periods $t$.
\item[$N$] Set of nodes $n$.
\item[$G$] Set of conventional generation units $g$.
\item[$F$] Set of stochastic power production units $f$.
\item[$L$] Set of loads $l$.
\item[$\Omega$] Set of scenarios $w$.
\item[$P$] Set of partitions $p$.
\item[$\Omega_p$] Set of scenarios $w$ in partition $p$.
\item[$F_n$] Set of stochastic power production units located at node $n$.
\item[$L_n$] Set of loads connected at node $n$.
\item[$G_n$] Set of conventional generation units located at node $n$.
\item[$M_n$] Set of nodes $m \in N$ that are connected to node $n$ by a transmission line.
\item[$\Omega^{\prime}_{p}$] Reduced set of scenarios $w$ in partition $p$.
\end{IEEEdescription}

\subsection{Parameters}
\begin{IEEEdescription}
\item[$C^{F}_{g},C^{V}_{g}$] Fixed/variable production cost of conventional generation unit $g$. 
\item[$C^{SU}_{g},C^{SD}_{g}$] Start-up/Shut-down cost of conventional generation unit $g$.
\item[$L_{l,t}$] Demand for load $l$ at time $t$.
\item[$RU_{g},RD_{g}$] Ramp-up/Ramp-down rate for conventional generation unit $g$.
\item[$UT_{g},DT_{g}$] Minimum-up/Minimum-down time for unit $g$
\item[$L^{UP}_{g},L^{DW}_{g}$] Number of time periods conventional generation unit $g$ must be online/offline counting from $t = 1$.
\item[$IS_{g}$] Initial status of unit $g$, equal to 1 if online at $t = 0$ and 0, otherwise.
\item[$ON_{g},OFF_{g}$] Number of time periods unit $g$ has been online/offline prior to $t = 1$.
\item[$X_{n,m}$] Reactance of line $n-m$.
\item[$F^{max}_{n,m}$] Maximum flow capacity of line $n-m$.
\item[$P^{max}_g,P^{min}_g$] Maximum/minimum power production of conventional generation unit $g$
\item[$P^{SU}_g,P^{SD}_g$]{ Maximum starting-up/shutting-down power production of conventional generation unit $g$.}
\item[$P^{IS}_g$] Power output of conventional unit $g$ at $t=0$.
\item[$C^{L}$] Cost of involuntary load curtailment.
\item[$W^{}_{f,t,w}$] Power production from stochastic generation unit $f$ at time $t$ in scenario $w$.
\item[$\pi_{w}$] Probability of scenario $w$.
\item[$\rho_{p}$] Weight associated with partition $p$.
\end{IEEEdescription}
%%FIRST STAGE VARIABLES
\subsection{First-stage variables}
\begin{IEEEdescription}
\item[$u_{g,t}$]{Binary variable equal to 1 if unit $g$ is online at time $t$ and 0, otherwise.}
\item[$y_{g,t}/z_{g,t}$]{Binary variable equal to 1 if  unit $g$ is starting up/shutting down at time  $t$ and 0, otherwise.}
\end{IEEEdescription}

%%SECOND STAGE VARIABLES
\subsection{Second-stage variables}
\begin{IEEEdescription}
\item[$P_{g,t,w}$]{Power produced by conventional generation unit $g$ in scenario $w$ at time $t$.}
\item[$L^{SH}_{l,t,w}$]{Power curtailment from load $l$ in scenario $w$ at time $t$.}
\item[$W^{SP}_{f,t,w}$]{Power curtailment from stochastic power production unit $f$ in scenario $w$ at time $t$.}
\item[$\delta_{n,t,w}$]{Voltage angle at node $n$, time $t$ and scenario~$w$.}
%\item{$F_{nq,t,w}$}{Flow from node n to node q in the balancing state at time $t$.}
\item[$\alpha$]{Auxiliary variable used in the scenario-based robust unit commitment formulation}
\item[$\theta_p$]{Auxiliary variable used in the hybrid unit commitment formulation}
\end{IEEEdescription}

\section{Introduction}
The increasing reliance on partly unpredictable renewable  power supply has prompted the revision of the procedures used for power system operations. This is the case, for example, of the tool used by system operators to decide the commitment of power plants, that is, to solve the so-called \emph{unit commitment problem} (UC). Two-stage stochastic programming \cite{birge2011introduction} and robust optimization \cite{ben2009robust} have become the most popular and explored techniques of optimization under uncertainty to improve unit-commitment decisions in terms of both cost-efficiency and system reliability. 

The formulation and solution of the unit commitment problem using either stochastic programming or robust optimization---the result of which is typically referred to as \emph{stochastic} and \emph{robust unit commitment}, respectively--- has been subject of numerous studies by the scientific community; see, for instance, \cite{takriti1996stochastic, bouffard2005market, wang2008security,morales2009economic, tuohy2009unit, wang2012chance, zhao2012robust, bertsimas2013adaptive, zeng2013solving}, among many others and variants.  

Essentially, the stochastic unit commitment problem (SUC) makes use of a probabilistic model for the uncertain input factors such as demand, equipment failures and partly-predictable renewable power production to minimize a certain quantile of the induced system cost distribution, usually, its expectation. Most often than not, this probabilistic model takes the form of a set of scenarios that describe plausible realizations of such random factors. In order for the stochastic solution to be reliable, the amount of scenarios that need to be considered must be large, which may render an intractable optimization problem, or carefully generated, which motivates the topic of \emph{scenario reduction techniques}  \cite{Heitsch, Morales2, Salva}. Furthermore, the probabilistic model from which these scenarios may be drawn may carry, in itself, some level of uncertainty as well. 

In contrast, the robust unit commitment problem (RUC) seeks a commitment plan that allows the system to withstand the worst-case realization of the uncertain factors at a minimum cost. While this approach saves the decision-maker from having to probabilistically characterize these factors, it may yield too conservative solutions as the worst-case scenario rarely occurs. 

In recent years, several methods have been proposed to make decisions under uncertainty that perform \emph{relatively well} under the premises of \emph{both} the stochastic and the robust approaches, that is, in expectation and for the worst case. Illustrative examples of these methods can be found in \cite{liu2003fuzzy,xu2009srccp, liubidding, fanzeres2015contracting}, where hybrid stochastic-robust solution strategies are developed for optimal air-quality and municipal solid-waste management, electricity trading for power microgrids and energy contracting for a portfolio of renewable power generation technologies, respectively. What makes all these solution strategies \emph{hybrid} is that some of the uncertain parameters are assumed to follow certain probability distributions, while others are solely known to belong to some uncertainty sets. 

Within the context of the unit commitment problem, we highlight the work in \cite{zhao2013unified} and \cite{zhao2015data}. More specifically, the authors in \cite{zhao2013unified} propose a mathematical formulation that delivers the unit-commitment plan that minimizes a user-controlled weighted sum of the expected and the worst-case costs.  The solution approach introduced in \cite{zhao2015data}, even if presented as a method to tackle the stochastic unit commitment problem, seeks to determine a unit-commitment plan that is robust against an \emph{ambiguous} probability distribution of renewable energy generation, an ambiguity that is the result of the always limited availability of data and that is modeled in practice as a vector of imperfectly known probabilities. Thus, as the amount of historical data increases, the ambiguity of such a probability distribution diminishes and so does the need for robustness and the degree of conservatism of the stochastic unit-commitment solution. This approach can also be regarded as a form of \emph{distributionally robust optimization} \cite{Gabrel2014471}.

Our work shares with \cite{zhao2015data} the aim of finding a solution to the stochastic unit commitment problem that is robust in some sense, but our motivation and the methodology we propose to this end are essentially different. We assume that the probability distributions of the uncertain parameters---in our case, the wind power production---are known, but that, as it normally occurs in practice, computationally tractability only allows us to solve the stochastic unit commitment problem for a scenario-based approximation of such distributions. In principle, we shall consider a large number of scenarios for this approximation to be accurate enough. In any case, we cluster these scenarios using the \emph{k-means clustering algorithm} \cite{wagstaff2001constrained}, which has been reported to feature good performance in similar contexts \cite{hobbs1999stochastic, watson2014new}. Each of the so-obtained clusters is referred to as a \emph{partition}. We then formulate and solve a two-stage unit-commitment problem that minimizes the expected value of the system operating costs, where the expectation is taken over the collection of worst-case scenarios within each partition. The probability assigned to each of these worst-case scenarios is equal to the probability of the partition it belongs to, which is, in turn, computed by summing up the probabilities of the scenarios that form part of the partition in question.

For convenience, we employ the term \emph{hybrid unit commitment problem} and the acronym HUC to refer to the proposed reformulation of the UC problem. This reformulation brings two major advantages, namely:

\begin{enumerate}
    \item It allows finding solutions to the two-stage unit commitment problem with different degrees of conservatism by changing the number of partitions. In fact, if only one partition is considered, the HUC delivers the robust unit-commitment solution. In contrast, if the number of partitions is made to coincide with the number of scenarios, the HUC solution boils down to the stochastic unit-commitment plan. 
    \item It is amenable to decomposition and paralellization at various levels and, hence, it can be efficiently solved. Indeed, each partition is processed in parallel. Furthermore, based on the column-and-constraint generation procedure described in~\cite{zeng2013solving}, we develop a master-subproblem decomposition scheme to find those scenarios within each partition that, considered together, result in the worst-case system operating cost for that partition. Finally, we solve the HUC problem for a substantially smaller set of scenarios, namely, those picked up from each partition.  In this sense, our solution approach works similarly to a scenario reduction technique that is robust against the  error intrinsic to the reduction process. 
\end{enumerate}

The remainder of this paper is organized as follows. Section \ref{SEC1} begins by providing mathematical formulations for the scenario-based two-stage stochastic and robust unit commitment problems, in that order, and finishes with the formulation of the proposed hybrid unit commitment problem. Furthermore, in this section we explain how we use the k-means clustering method to construct the partitions in the HUC model and how these can be employed to control the degree of conservatism of the resulting unit-commitment plan. Section \ref{BDA} introduces the proposed parallelization and decomposition strategy to solve the HUC problem. Section \ref{CaseStudies} analyzes and discusses results from two case studies based on standard IEEE power systems. Finally, in Section \ref{Conclusion} the main conclusions of our work are summarized, including possible avenues for future research.

\section{Mathematical Formulation}
\label{SEC1}
In the two-stage unit commitment problem under uncertainty, decision variables are divided into two groups. The first group constitutes the commitment plan itself and consists of the 0/1 variables $u_{g,t}$, $y_{g,t}$, $z_{g,t}$, which determine the on/off status, the start-up, and the shutdown of generating unit $g$ in time period $t$, respectively. These decisions are to be made, in general, one day in advance of the actual delivery of electricity and, in any case, before the realization of the uncertain factors. In this paper, we consider for simplicity that the system uncertainty stems only from the wind power production, which is modeled as a finite set $\Omega$ of scenarios $W^{}_{f,t,w}$ with $w \in \Omega$.  

\indent The second-stage decision variables, namely, $P_{q,t,w}$, $L^{SH}_{l,t,w}$, $W^{SP}_{f,t,w}$ and $\delta_{n,t,w}$ determine the economic dispatch of the conventional generating units, the amount of load that is involuntarily shed, the amount of wind power production that is curtailed, and the voltage angles at the network nodes, respectively. These variables adapt to the specific realization of the uncertainty and as such, are augmented with the scenario index $w$. 
\\  
\indent We start by providing the mathematical  formulation of the two-stage \emph{stochastic} unit commitment problem. In all cases, we consider that the marginal production cost of the wind generation is zero.
\subsection{Two-Stage Stochastic Unit Commitment (SUC)}
\label{SUC}
The two-stage stochastic unit commitment problem can be formulated as follows:
\begin{align}
\underset{\mathcal{H}, \mathcal{W}}{\text{minimize}}\hspace{0.2cm}
  &\sum_{t\in T}\sum_{g\in G}\left(C^{F}_{g} u_{g,t}+C^{SU}_{g} y_{g,t}+C^{SD}_{g} z_{g,t}\right) \label{9} \\ \notag & + \sum_{w\in \Omega}\pi_{w}\left[\sum_{t\in T}\sum_{g\in G}C^{V}_{g} P_{g,t,w}+\sum_{t\in T}\sum_{l\in L} C^{L} L^{SH}_{l,t,w}\right]  \allowdisplaybreaks \\[0.2cm]
\text{s.t.}\hspace{0.2cm}
& y_{g,t}-z_{g,t} = u_{g,t}-u_{g,t-1} \label{10} \\ \hspace{1cm} & (\forall{g},\forall{t} \in \{2,...,T\}) \notag  
  \allowdisplaybreaks \\[0.2cm]
& y_{g,t}-z_{g,t} = u_{g,t}-IS_{g} \label{11} \\ \hspace{1cm} & (\forall{g},\forall{t} \in \{1\})  \notag
  \allowdisplaybreaks \\[0.2cm]
& y_{g,t}+z_{g,t} \leq 1   \label{12} \\ \hspace{1cm} & (\forall{g},\forall{t} \in \{1,...,T\})  \notag
  \allowdisplaybreaks \\[0.2cm]
& u_{g,t}=IS_{g}  \label{13}  \\ \hspace{1cm} & (L^{UP}_{g}+L^{DW}_{g} > 0, \forall g, \forall t \leq L^{UP}_{g}+L^{DW}_{g}) \notag
\allowdisplaybreaks \\[0.2cm]
& \sum^{t}_{\tau=t-UT_{g}+1} y_{g,\tau} \leq u_{g,t}  \label{14}   \\ \hspace{1cm} & (\forall g, \forall t > L^{UP}_{g}+L^{DW}_{g})\notag
\allowdisplaybreaks \\[0.2cm]
& \sum^{t}_{\tau=t-DT_{g}+1} z_{g,\tau} \leq 1-u_{g,t}  \label{15} \\ \hspace{1cm} & (\forall g, \forall t > L^{UP}_{g}+L^{DW}_{g})\notag
\allowdisplaybreaks \\[0.2cm]
& \sum_{g\in G_n}P_{g,t,w}-\sum_{l\in L_n}L_{l,t}+\sum_{l\in L_n}L^{SH}_{l,t,w}+\sum_{f\in F_n}W^{}_{f,t,w}  \label{16} \\[-0.2cm] \notag & -\sum_{f\in F_n}W^{SP}_{f,t,w}=\sum_{q\in Q_n}\frac{(\delta_{n,t,w}-\delta_{q,t,w})}{X_{n,m}}  \\ \hspace{1cm} & (\forall{n},\forall{t},\forall{w}\in \Omega) \notag &
 \allowdisplaybreaks \\[0.2cm]
% & \sum_{g\in G_n}P_{g,t,w}-\sum_{l\in L_n}L_{l,t}+\sum_{l\in L_n}L^{SH}_{l,t,w}+\sum_{f\in F_n}W^{}_{f,t,w}  \label{16} \\ \notag & -\sum_{f\in F_n}W^{SP}_{f,t,w}=\sum_{q\in Q_n}F_{nq,t,w}  \\ \hspace{1cm} & (\forall{n},\forall{t},\forall{w}\in \Omega) \notag &
% \allowdisplaybreaks \\[0.2cm]
 & P_{g,t,w} \leq P^{max}_{g}u_{g,t}  \label{17}  \\ \hspace{1cm} & (\forall{g},\forall{t} ,\forall{w}\in \Omega) \notag &
 \allowdisplaybreaks \\[0.2cm]
& P_{g,t,w} \geq P^{min}_{g}u_{g,t}   \label{18}  \\ \hspace{1cm} & (\forall{g},\forall{t} ,\forall{w}\in \Omega)  \notag &
 \allowdisplaybreaks \\[0.2cm]
& P_{g,t,w} \leq (P^{IS}_{g}+RU_{g})u_{g,t}  \label{19} \\ \hspace{1cm} & (\forall{g},\forall{t} \in \{1\},\forall{w}\in \Omega) \notag &
  \allowdisplaybreaks \\[0.2cm]
  &   P_{g,t,w} \geq (P^{IS}_{g}-RD_{g})u_{g,t}  \label{20} \\ \hspace{1cm} & (\forall{g},\forall{t} \in \{1\},\forall{w}\in \Omega) \notag &
  \allowdisplaybreaks \\[0.2cm]
& P_{g,t,w}-P_{g,t-1,w}\leq  (2-u_{g,t-1}-u_{g,t})P^{SU}_{g}  \label{21} \\ \notag & +(1+u_{g,t-1}-u_{g,t})RU_{g} \\ \hspace{1cm} & (\forall{g},\forall{t} \in \{2,...,T\},\forall{w}\in \Omega) \notag &
 \allowdisplaybreaks \\[0.2cm]
& P_{g,t-1,w}-P_{g,t,w}\leq (2-u_{g,t-1}-u_{g,t})P^{SD}_{g}  \label{22} \\ \notag & +(1-u_{g,t-1}+u_{g,t})RD_{g} \\  \hspace{1cm} & (\forall{g},\forall{t} \in \{2,...,T\},\forall{w}\in \Omega) \notag &
 \allowdisplaybreaks \\[0.2cm]
& L^{SH}_{l,t,w} \leq L_{l,t}  \label{23}  \\ \hspace{1cm} & (\forall{l},\forall{t} ,\forall{w}\in \Omega) \notag &
 \allowdisplaybreaks \\[0.2cm]
& W^{SP}_{f,t,w} \leq W^{}_{f,t,w}  \label{24}  \\ \hspace{1cm} & (\forall{f},\forall{t} ,\forall{w}\in \Omega) \notag &
 \allowdisplaybreaks \\[0.2cm]
& -F^{max}_{n,m} \leq \frac{(\delta_{n,t,w}-\delta_{q,t,w})}{X_{n,m}}\leq F^{max}_{n,m}   \label{25} \\ \hspace{1cm} & (\forall{n,q\in Q_n},\forall{t} ,\forall{w}\in \Omega) \notag &
 \allowdisplaybreaks \\[0.2cm]
& P_{g,t,w}, L^{SH}_{l,t,w}, W^{SP}_{f,t,w} \geq 0  \label{26}  \\ & (\forall{g},\forall{l},\forall{f},\forall{t},\forall{w}\in \Omega) \notag
 \allowdisplaybreaks \\[0.2cm]
& u_{g,t},y_{g,t},z_{g,t} \in \{0,1\}  \label{27}  \\ & \notag (\forall{g},\forall{t})
\end{align}
where $\mathcal{H} = \left\{u_{g,t}, y_{g,t}, z_{g,t}\right\}$ and $\mathcal{W} = \left\{P_{g,t,w}, L^{SH}_{l,t,w}, W^{SP}_{f,t,w}, \delta_{n,t,w}:\omega \in \Omega \right\}$ are the sets of here-and-now and wait-and-see decisions, respectively. Furthermore, following \cite{carrion2006computationally}, the initial state conditions are given by 
\begin{align*}
IS_{g}=\left\lbrace\begin{array}{ll}
 1 & \textup{if } ON_{g}>0\\
 0 & \textup{if } ON_{g}=0
\end{array}
\right.
\end{align*}
\begin{align*}
L^{UP}_{g}&=\text{min}\{T,(UP_{g}-ON_{g})IS_{g}\} \\[0.2cm]
L^{DW}_{g}&=\text{min}\{T,(DT_{g}-OFF_{g})(1-IS_{g})\}
\end{align*}

Problem~\eqref{9}--\eqref{27} takes the form of a standard two-stage unit commitment formulation, which is similar, to a large extent, to those provided in the numerous works on the topic, see, for instance,~\cite{morales2009economic, Papavasiliou} and references therein. The objective is to minimize the expected system operating cost~\eqref{9}, which is made up of the no-load, start-up, shutdown, and variable production costs of the conventional generating units, and the cost of involuntarily load curtailment, in that order. Equations~\eqref{10}--\eqref{12} model the changes in the on/off-commitment status of the power plants as these are started up or shutdown throughout the scheduling horizon, while \eqref{13}--\eqref{15} impose their minimum up- and down-time requirements. Equalities~\eqref{16} constitute the set of nodal power balance equations according to a DC power flow model. The maximum and minimum power outputs of the generating units are enforced by \eqref{17} and \eqref{18}, respectively, and their ramping limits through \eqref{19}--\eqref{22}, as in \cite{wang2012chance} and \cite{zhao2012robust}. The sets of inequalities \eqref{23} and \eqref{24} limit the involuntary load curtailment and the wind power spillage to the eventual power that is consumed and the eventual wind power that is produced, respectively. The set of equations~\eqref{25} guarantee compliance with the transmission capacity limits. Finally, constraints~\eqref{26} and \eqref{27} constitute variable declarations.

%\cite[Ch. 7]{morales2013integrating}, except for the fact that we treat the power dispatch of conventional generating units as recourse decisions.  In particular The objective function, to be minimized 

%However, in our model we present some differences regarding the first-stage decision variables. In \cite{morales2013integrating} the power production committed is included as first-stage variables and therefore upward and downward regulation has to be added, whilst in our model, only the commitment status is included in the first stage variables and therefore, there is no need for regulation since the power production will be dispatched afterwards. More differences are the ramping constraints (\ref{21})-(\ref{22}) formulation, based on  and the initial state conditions, . 

\subsection{Two-Stage Robust Unit Commitment (RUC)}
\label{RUC}
In this section we present a scenario-based formulation of the two-stage robust unit commitment problem, which writes as follows:
\begin{align}
&\underset{\mathcal{H}, \mathcal{W}, \alpha}{\text{minimize}}\hspace{0.2cm}
  \sum_{t\in T}\sum_{g\in G}\left(C^{F}_{g} u_{g,t}+C^{SU}_{g} y_{g,t}+C^{SD}_{g} z_{g,t}\right) + \alpha \label{37} \\
& \text{s.t.}\hspace{0.2cm}  \alpha \geq \sum_{t\in T}\sum_{g\in G}C^{V}_{g} P_{g,t,w}+\sum_{t\in T}\sum_{l\in L} C^{L} L^{SH}_{l,t,w}, \enskip \forall{w} \in \Omega \label{38} \\
& (\ref{10})-(\ref{27}) \label{75}
\end{align}  
In (\ref{37})--(\ref{75}), the auxiliary variable $\alpha$ equals the worst-case dispatch cost at the optimum. Note that this variable is bounded from below by a finite set of linear constraints (\ref{38}), one per scenario, that involve the second-stage decision variables $P_{g,t,w}$ and $L^{SH}_{l,t,w}$. Thus, the objective of the scenario-based RUC problem (\ref{37})--(\ref{75}) is to minimize the total system operating cost for the worst-case scenario of the uncertainty.

In the following section, we introduce the proposed hybrid formulation of the two-stage unit commitment problem under uncertainty.

\subsection{Hybrid Unit Commitment Problem (HUC)}
\label{HUC}

The two-stage unit commitment formulation that we propose is based on splitting the finite set of scenarios $\Omega$ into $k$ partitions, each with a probability equal to the sum of probabilities of the scenarios that form part of it. The objective is then to minimize the expected system operating cost over the scenarios that deliver the worst-case dispatch cost within each partition. These scenarios, one per partition, are assigned a probability of occurrence equal to the probability of the partition they belong to. 

Let $\Omega=\{1,...,\lambda\}$ denote the original scenario set, where $\lambda$ is the total number of scenarios. These scenarios are then clustered into $k$ partitions with $P=\{1,...,k\}$ being the partition set. For ease of notation, we define the series of subsets $\Omega_{1}, \ldots, \Omega_{p}, \ldots, \Omega_{k}$, with $\Omega_{i} \bigcap \Omega_{j} = \emptyset$ for all $i \neq j$ and $\Omega_1 \bigcup \ldots \Omega_p \bigcup \ldots \Omega_k = \Omega$, such that $\Omega_p$ is comprised of all the scenarios $w \in \Omega$ that belong to partition $p \in P$.

Since $\Omega$ is a discrete probability space with probability measure $\pi_{w}$, $w = 1, \ldots, \lambda$, the probability $\rho_{p}$ associated with each partition $\Omega_p$, $p \in P$, depends on the number of scenarios that pertain to it and is computed as: 
\begin{align}
\rho_{p}=\sum_{w \in \Omega_p} \pi_w & \quad \forall p\in P \label{39}  \allowdisplaybreaks  \\
\sum_{p \in P} \rho_p =1
\end{align}

Therefore, the proposed hybrid two-stage unit commitment problem writes as follows:

\begin{align}
&\underset{\mathcal{H}, \mathcal{W},\theta_p}{\text{min.}}\hspace{0.2cm}
  \sum_{t\in T}\sum_{g\in G}\left(C^{F}_{g} u_{g,t}+C^{SU}_{g} y_{g,t}+C^{SD}_{g} z_{g,t}\right) + \allowdisplaybreaks  \notag \\[0.3cm]
& \hspace{0.8cm} + \sum_{p \in P} \rho_p \theta_p \label{46}  \allowdisplaybreaks  \\[0.3cm]
& \text{s.t.}\hspace{0.2cm}  \theta_p \geq \sum_{t\in T}\sum_{g\in G}C^{V}_{g} P_{g,t,w}+\sum_{t\in T}\sum_{l\in L} C^{L} L^{SH}_{l,t,w}, \label{47}   \\  & \quad\enskip (\forall p \in P, \hspace{0.1cm} \forall w \in \Omega_{p}) \notag  \allowdisplaybreaks  \\[0.3cm]
 & \quad\enskip (\ref{10})-(\ref{27}) \label{76}
\end{align}

The auxiliary variable $\theta_p$, one per partition $p \in P$, equals the worst-case dispatch cost within partition $p$, in a similar way as the auxiliary variable $\alpha$ does in the robust unit commitment formulation (\ref{37})--(\ref{75}) for the whole set of scenarios $\Omega$. This way, problem (\ref{46})--(\ref{76}) is expected to yield a unit-commitment plan that is ``in between'' the robust and the stochastic unit-commitment solutions in terms of the expected and the worst-case system operating cost. Furthermore, the closeness of the HUC solution to the stochastic and robust unit-commitment plans, and consequently its degree of conservatism, are controlled by the number $k$ of partitions or clusters into which the scenarios are grouped. Indeed, if the number of partitions equals the number of scenarios, that is, $k = \lambda$, the HUC model (\ref{46})--(\ref{76}) reduces to~\eqref{9}--\eqref{27} and the stochastic solution is obtained. In contrast, if only one single partition is considered ($k = 1$), we have that $\Omega_1 = \Omega$ and therefore, problem (\ref{46})--(\ref{76}) boils down to the scenario-based robust unit-commitment formulation (\ref{37})--(\ref{75}). As a result, the HUC solution coincides with the robust solution in such a case.

Hence, we can increase  the degree of conservatism of the HUC solution by diminishing the number of partitions, and vice versa. For $1 < k < \lambda$, however, how efficiently the HUC solution transits from the robust to the stochastic unit-commitment plan, as $k$ increases, depends on the performance of the clustering technique. We use the k-means algorithm \cite{wagstaff2001constrained}, which is an efficient non-hierarchical method to cluster a data set into $k$ groups. In our case, the k-means algorithm assigns each scenario $w \in \Omega$ to the partition $\Omega_p$, with $p\in P$, with the nearest mean. The k-means algorithm has been reported to showcase the best performance in a probabilistic production cost model in \cite{hobbs1999stochastic} and in transmission and generation expansion planning in \cite{watson2014new}, compared to several other clustering techniques. Furthermore, this previously reported evidence agrees with what we observe in our numerical experiments.

\section{Solution Strategy: Parallelization and Decomposition}

It is well known that the unit commitment problem is mixed-integer, NP-hard, and generally requires long solution times. This is especially true for realistic instances of the two-stage unit commitment problem under uncertainty.
In the following we describe the parallelization-and-decomposition scheme that we have designed to efficiently solve the proposed HUC formulation (\ref{46})--(\ref{76}). For ease of exposition, we divide this description in two parts. In the first one, we explain how problem (\ref{46})--(\ref{76}) is decomposed per partition and scenario, while in the second part we elaborate on how the solution to the decomposed problem is parallelized.

\label{BDA}
\subsection{Problem Decomposition via Column-and-constraint Generation}

Let us consider a certain partition $p \in P$ that comprises the subset of scenarios $\Omega_p$. Note that, for determining the optimal solution to the HUC problem (\ref{46})--(\ref{76}), we only need those (hopefully few) scenarios in $\Omega_p$ that deliver the worst-case dispatch cost within partition $p$ for any feasible unit-commitment plan. Let $\Omega'_{p} \subset \Omega_p$ denote the subset of those scenarios. The idea is to find such a reduced set $\Omega'_{p}$, a task that can be performed in parallel for each partition $p \in P$.

To build $\Omega'_{p}$ from $\Omega_p$, the latter being the outcome of the k-means clustering algorithm, we develop a master-subproblem decomposition scheme based on the column-and-constraint generation procedure described in \cite{zeng2013solving}. In the sequel we will refer to this decomposition scheme as \emph{Primal Cut Algorithm} after the solution strategy introduced in \cite{zhao2012robust} whereby the master problem is gradually enlarged with the addition of cuts expressed in terms of the primal variables.

 The master problem is a mixed-integer programming problem that involves both first-stage and second-stage decision variables and that has the following form at iteration $i$ of the column-and-constraint generation algorithm:
\begin{align}
&\underset{\mathcal{H}^{i}, \mathcal{W}^{i},\theta_p}{\text{minimize}}\hspace{0.2cm}
  \sum_{t\in T}\sum_{g\in G}\left(C^{F}_{g} u_{g,t}^{i}+C^{SU}_{g} y_{g,t}^{i}+C^{SD}_{g} z_{g,t}^{i}\right) + \theta_p \label{48}  \allowdisplaybreaks  \\[0.3cm]
& \text{s.t.}\hspace{0.2cm} \quad\enskip (\ref{10})-\eqref{15},(\ref{27})\\[0.3cm]  
 & \theta_p \geq \sum_{t\in T}\sum_{g\in G}C^{V}_{g} P_{g,t,w}^{i}+\sum_{t\in T}\sum_{l\in L} C^{L} L^{sh,i}_{l,t,w}, \quad \forall w \in \Omega^{\prime i}_{p} \label{49}\\[0.3cm]
 & \quad\enskip \eqref{17}-\eqref{26}, \quad\forall w \in \Omega^{\prime i}_{p} \label{50}
\end{align}
where $\mathcal{H}^{i} = \left\{u_{g,t}^{i},y_{g,t}^{i},z_{g,t}^{i}\right\}$ and $\mathcal{W}^{i} = \big\{P_{g,t,w}^{i}, L^{sh,i}_{l,t,w},W^{sp,i}_{f,t,w},\delta_{n,t,w}^{i}: \omega \in \Omega_{p}^{'i}\big\}$. Note that $\Omega_{p}^{'0} = \emptyset$. As the algorithm proceeds, $\Omega_{p}^{'i}$ is augmented with those possibly few scenarios $\omega \in \Omega_p$ that are needed to reconstruct the partition-worst-case recourse cost as a function of the first-stage decision variables $u_{g,t}^{i}$, $y_{g,t}^{i}$, and $z_{g,t}^{i}$ in the form of \eqref{49}--\eqref{50}. 

Constraint \eqref{49} can be interpreted as a primal cut, as compared to those cuts that are constructed from dual information, as it is the case, for example, of a standard Benders cut.

The subproblems are linear programming problems (LP) that determine the second-stage decision variables $P_{g,t,w}^{i}$, $L^{sh,i}_{l,t,w}$, $W^{sp,i}_{f,t,w}$, and $\delta_{n,t,w}^{i}$ with $u_{g,t}^{i}$, $y_{g,t}^{i}$, and $z_{g,t}^{i}$ fixed at the values given by the master problem. A subproblem in the form of \eqref{51}--\eqref{52}
is solved for each scenario $\omega \in \Omega_p$.
\begin{align}
\underset{\mathcal{W}^i_w }{\text{minimize}}\hspace{0.2cm}
&  \sum_{t\in T}\sum_{g\in G}C^{V}_{g} P^{i}_{g,t,w}+\sum_{t\in T}\sum_{l\in L} C^{L} L^{sh,i}_{l,t,w} \label{51}\\[0.3cm]
 \text{s.t.}\hspace{0.2cm} 
  &(\ref{16})-(\ref{26}) \label{52}
\end{align} 
where $\mathcal{W}^i_w = \left\{P_{q,t,w}^{i}, L^{sh,i}_{l,t,w}, W^{sp,i}_{f,t,w}, \delta_{n,t,w}^{i}\right\}$.

The scenario $w'$ for which the associated subproblem~\eqref{51}--\eqref{52} yields the highest dispatch cost or is infeasible is used to construct a set of primal constraints in the form of \eqref{49}--\eqref{50} that is added to the master problem by setting $\Omega_{p}^{'i+1} = \Omega_{p}^{'i} \bigcup \{w'\}$.
It is worth noticing, however, that subproblem infeasibility is not a concern in our case due to the possibility of shedding load and spilling wind.

%In order to avoid rework every time that the PCA must solve the master problem, we implement the proposed master-subproblem decomposition scheme using callbacks. This means that we do not solve the master problem to optimality, but instead we pass the subproblems the first integer feasible solution (the firs incumbent solution) that the branch-and-cut algorithm used to solve the master problem finds. This remarkably speeds up the solution process. 

One instance of the primal cut algorithm is run for each partition $p \in P$ in parallel. Each of these instances works, therefore, with one master problem and a number of subproblems equal to the number of scenarios in each partition, that is, equal to $card\left(\Omega_{p}\right)$. Furthermore, each instance of the algorithm concludes by delivering the set of selected scenarios  $\Omega_{p}^{\prime} \subset \Omega_p$ for partition $p$. The last step of our solution strategy consists then in solving the HUC problem (\ref{46})--(\ref{76}) where $\Omega_{p}$ is replaced with the reduced scenario set $\Omega^{\prime}_{p}$.  

We describe below how our solution strategy proceeds step by step. 

\begin{enumerate}
  \item Choose the number $k$ of partitions and apply the k-means clustering method to the complete set of scenarios $\Omega$ in order to assign each scenario to a certain partition $p$. 
  \vspace{0.1cm} 
      \item Create one instance of the primal cut algorithm for each partition $p \in P$.
        \vspace{0.1cm} 
  \item Initialization: Set $i = 0$ and $\Omega_{p}^{'0} = \emptyset$. 
   \vspace{0.1cm} 
  \item Solve the master problem (MP). Return the optimal solution found by the branch-and-cut algorithm and denote this solution by $(u_{g,t}^{i}, y_{g,t}^{i}, z_{g,t}^{i})$. Calculate a lower bound LB as $\sum_{t\in T}\sum_{g\in G}(C^{F}_{g} u_{g,t}^{i} + C^{SU}_{g} y_{g,t}^{i} + C^{SD}_{g} z_{g,t}^{i})+ \theta_p$.
    \vspace{0.1cm} 
  \item Solve the subproblems (SP) with the first-stage decision variables fixed at $(u_{g,t}^{i}, y_{g,t}^{i}, z_{g,t}^{i})$. Once the SP are solved, the scenario $w'$ associated with the subproblem that yields the highest dispatch cost is identified and included into the reduced set $\Omega^{\prime i+1}_{p}$, i.e., $\Omega^{\prime i+1}_{p} = \Omega^{\prime i}_{p} \bigcup \{w'\}$. Compute an upper bound $UB$ as  $\sum_{t\in T}\sum_{g\in G}\left(C^{F}_{g} u_{g,t}^{i}+C^{SU}_{g} y_{g,t}^{i}+C^{SD}_{g} z_{g,t}^{i}\right) +  \sum_{t\in T}\sum_{g\in G}C^{V}_{g} P_{g,t,w'}^{i}+\sum_{t\in T}\sum_{l\in L} C^{L} L^{sh,i}_{l,t,w'}$.   
 \vspace{0.1cm} 
    \item Convergence check: If $\mid UB-LB\mid \leq \epsilon$, being $\epsilon$ a user-specified tolerance value, the iterative process stops. If $\mid UB-LB\mid > \epsilon$, then set $i := i+1$ and go to step 4.
  \vspace{0.1cm} 
      \item Once all the instances of the primal cut algorithm have converged, the HUC problem (\ref{46})--(\ref{76}) is solved  for all $p \in P$ and for all $w \in \Omega^{\prime}_{p}$. The reduced set $\Omega^{\prime}_{p}$ is made up of those scenarios $w \in \Omega_{p}$ that determine the worst-case dispatch cost within partition $p$.     
\end{enumerate}
A pseudocode for the proposed decomposition scheme, which hereinafter we refer to as \emph{Scenario Partition and Decomposition Algorithm} (SPDA), is provided in Algorithm \ref{CHalgorithm}. For ease of notation, let $x$ ($x^{i}$) denote the vector of first-stage variables (at iteration $i$).
\begin{algorithm}
\caption{Scenario Partition and Decomposition Algorithm (\textbf{SPDA})}
\label{CHalgorithm}
\begin{algorithmic}[1]
\State Choose $k$ and apply k-means to $\Omega$. 
\For{all $p$ \Pisymbol{psy}{206} $P$ }
 \State \textbf{Set} $i:= 0$ and $\Omega_{p}^{'0} = \emptyset$
\Repeat
\State \textbf{Solve} MP  $\forall w \in \Omega_{p}^{\prime i}$
\State Return optimal solution $x^{i}$
\State Compute $LB$
\State \textbf{Set} $x := x^{i}$ and \textbf{solve} SP $\forall w \in \Omega_{p}$ 
\State Compute UB
\State Identify worst-case scenario $w'$ 
\State \textbf{Set} $\Omega_{p}^{\prime i+1}:= \Omega_{p}^{\prime i} \bigcup  \{w'\}$
\State \textbf{Set} $i:= i +1$
\Until $\mid UB-LB\mid \leq \epsilon$ 
\State \textbf{Set} $\Omega_{p}^{\prime}:= \Omega_{p}^{\prime i-1}$
\EndFor
\State \textbf{Solve} HUC $\forall p \in P, \hspace{0.1cm} \forall w \in \Omega_{p}^{\prime}$ 
\end{algorithmic}
\end{algorithm}
 
Notice that SPDA works as a scenario reduction technique that retains the most detrimental scenarios in terms of system operating cost. This confers robustness to the solution of the proposed HUC problem. Moreover, the last command line in SPDA, which involves solving the HUC model for the reduced scenarios sets $\Omega_{p}^{\prime}, \forall p \in P$, could be carried out as well via further decomposition (see, for instance, \cite{Papavasiliou2015}), although this possibility has not been explored in this paper. 

\subsection{Parallelization of the Solution Algorithm}

In SPDA both the outer ``for-loop''and the solution to the w-indexed subproblems are amenable to parallelization. For this purpose, we make use of the DTU High Performance Computing (HPC) Facility \cite{Centr57:online}. We create $k$ jobs, each representing an instance of the primal cut algorithm for each of the $k$ partitions into which we divide the scenario set $\Omega$. These jobs are simultaneously submitted to the HPC Cluster, where they are concurrently executed, as there is no need for communication in between the workers (nodes or cores). %%The HPC Cluster manages and distributes the available resources to the workers in an opportunistic manner according to their workload and the required memory. 

%%Besides, each of these $k$ jobs involves solving a different number of subproblems for each iteration of the running instance of the primal cut algorithm. The solution processes of these subproblems are also handled in parallel by the workers (based on their workload and memory requirements as well). To generate the subproblems and execute them in a distributed fashion within the HPC Cluster, we employ the GAMS Grid Facility \cite{bussieck2009grid}. 

We submit each of the $k$ jobs to a different node, using the same amount of resources per node. Within every node, the subproblems are solved in a multi-threaded environment using the Gather-Update-Solve-Scatter Facility in GAMS \cite{GAMSD96:online}. This tool allows treating each subproblem, one per scenario in the partition under consideration, as a different parametrization of the same linear programming model, which is then generated only once by GAMS. Likewise, the solutions to all the subproblems (or portions thereof) are retrieved back to GAMS in a single transaction.

%although the subproblems can be submitted in parallel through different cores, it turns out more effective to solve them in a multi-threading environment using the Gather-Update-Solve-Scatter Facility in GAMS \cite{GAMSD96:online}. This tool allows to simulate each of the scenarios that compose the partition under consideration by generating one single subproblem and retrieve portions of the solutions back in the model environment sequentially without returning to GAMS. 
%whereby only one single subproblem model is generated and then run and solved for each of the scenarios that compose the partition under consideration.  

\begin{figure}[h]
\centering
    \includegraphics[scale=0.8]{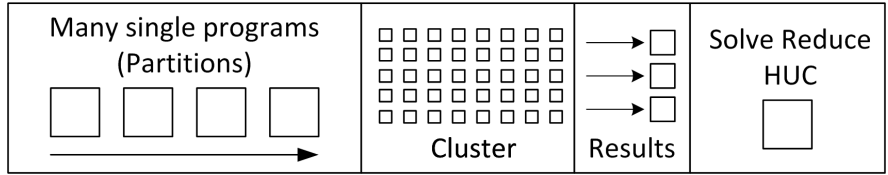}
     \caption{Representation of a cluster that runs many core independent jobs simultaneously describing the SPDA parallelization process.}
     \label{Paralelization}
    \end{figure}

Figure \ref{Paralelization} provides a graphical illustration of how SPDA is parallelized and run using a multi-machine configuration within the DTU HPC Cluster.  
%
%The resolution of the PBDA algorithm is an scenario independent process. Consequently, all the PBDA instances and subproblems can be run concurrently, in the sense that it is not necessary a communication process between the workers (nodes or cores). 
%The purpose is to allow the cluster to act as a computer farm. 
%The job is not run in a local workstation but it is submitted to the cluster for the execution, afterwards, the cluster will manage the resources needed to solve the job in an opportunistic way according to the memory required. Thus, we define one single program for each partition and we send them simultaneously to the cluster. Once in the cluster, each program performs a different number of subproblems that are solved in parallel by the workers. To generate the subproblems and execute them in a parallel distribution we use 

\section{Case Studies}\label{CaseStudies}
In the following the quality of the unit-commitment plan provided by the proposed HUC formulation is first tested on the IEEE 14-node power system \cite{Power25:online} and the IEEE 3-Area RTS-96 system \cite{TheRE28:online} in Section~\ref{ResultsQuality}. The latter power system is then used to evaluate the performance of SPDA in Section~\ref{ResultsPerformance}. 

The IEEE 14-node system comprises 14 nodes, 5 generators, 20 lines and 11 loads. We also add one wind farm to node 5, whose power production is modeled by the ten scenarios provided in Table \ref{WPP}.

The IEEE 3-Area RTS-96 system consists of 72 nodes, 96 generators, 107 lines, and 51 loads. Besides, we add 15 wind farms of 200-MW capacity each and location given by Table \ref{wfloc}. Thus, the wind power capacity represents 29\% of the total generating capacity installed in the power system. The technical characteristics of the generating units, the demand and the transmission lines are available online \cite{TheRE28:online}. The wind power production scenarios used for this case study come from \cite{pinson2013wind}, where the spatio-temporal dependencies of wind power generation are considered. More specifically, the study in \cite{pinson2013wind} provides 100 scenarios of wind power production that were generated for 15 control areas and 43 lead times in western Denmark. However, for this work, only 50 equiprobable scenarios and 24 lead times are considered.

We set the MIP tolerance gap to $0$ in all the simulations pertaining to the IEEE 14-node system, while we allow for a MIP tolerance gap of up to $10^{-3}$ $(0.1\%)$ in those numerical experiments carried out on the IEEE 3-Area RTS-96 system.

\subsection{Results Obtained and Effect of Partitions}\label{ResultsQuality}
 We first consider the IEEE 14-node system and compare the unit-commitment plans resulting from the proposed HUC formulation with those obtained using the method proposed in \cite{zhao2013unified}, which we refer to as \emph{Zhao model} hereafter. The comparison is conducted for a number of partitions in the HUC problem ranging from 1 to 10 and a number of values for the weighting factor in Zhao model varying from 0 to 1. This process prompts two different unit-commitment plans only, namely, UCP1 and UCP2, whose quality is measured in terms of both the expected and the worst-case system operating cost and indicated in Table  \ref{table:zhaopartitions}. Recall that the pure stochastic and robust unit-commitment solutions are obtained for a number of partitions equal to 10 and 1 in the proposed HUC formulation, and for values of the weighting factor in Zhao model equal to 0 and 1, respectively. 
\\

 \begin{table}[h]
\begin{center}
\caption{Wind power production scenarios [MW] for the IEEE 14-node system.}
 \scalebox{1}[1]{
\begin{tabular}{|c|c|c|c|c|c|c|c|c|c|c|}
\hline
  & w1 & w2 & w3 &w4 & w5 & w6 & w7 & w8 & w9 & w10 \\
\hline\hline
t1	&	66.45	&	67.06	&	65.39	&	57.12	&	42.9	&	37.8	&	38.64	&	29.12	&	27.96	&	22.35 \\
\hline
t2	&	107.4	&	102.34	&	86.71	&	85.2	&	70.29	&	65.79	&	59.2	&	51.52	&	40.26	&	35.6 \\
\hline
t3	&	113.25	&	103.32	&	100.23	&	90.6	&	74.47	&	66.24	&	63.76	&	51.45	&	47.82	&	37.25 \\
\hline
t4	&	127.5	&	102.2	&	113.1	&	100.56	&	90.2	&	78.57	&	61.04	&	58.03	&	46.5	&	40.85 \\
\hline
t5	&	124.2	&	113.68	&	106.08	&	109.2	&	94.6	&	78.48	&	68.32	&	55.3	&	51.84	&	37.95 \\
\hline
t6	&	128.25	&	113.4	&	110.63	&	88.92	&	93.72	&	78.57	&	68.8	&	55.86	&	50.22	&	41.75 \\
\hline
t7	&	152.1	&	147.98	&	127.14	&	116.04	&	101.31	&	91.71	&	78.48	&	73.85	&	60.78	&	52.65 \\
\hline
t8	&	158.55	&	139.16	&	117.65	&	120.36	&	109.89	&	80.82	&	73.68	&	69.09	&	63.48	&	53.7 \\
\hline
t9	&	120.15	&	119.14	&	107.77	&	93.24	&	81.4	&	70.65	&	58.64	&	51.24	&	46.2	&	42.5 \\
\hline
t10	&	85.8	&	86.94	&	81.9	&	81.6	&	72.93	&	60.03	&	47.44	&	43.05	&	37.98	&	32.55 \\
\hline
t11	&	148.2	&	141.68	&	133.51	&	128.88	&	105.05	&	89.37	&	82.16	&	67.41	&	52.74	&	52.7 \\
\hline
t12	&	131.25	&	129.08	&	113.75	&	104.04	&	98.01	&	86.94	&	69.6	&	68.46	&	61.38	&	45.6 \\
\hline
t13	&	135.45	&	120.12	&	117.65	&	102.36	&	93.61	&	71.28	&	62.4	&	59.5	&	46.98	&	39.55 \\
\hline
t14	&	118.5	&	109.9	&	112.45	&	96	&	81.62	&	76.14	&	59.36	&	52.36	&	51.48	&	41.2 \\ 
\hline
t15	&	110.25	&	114.52	&	104.39	&	90.84	&	92.51	&	72.18	&	62.4	&	50.96	&	47.46	&	39.9 \\
\hline
t16	&	44.4	&	43.4	&	41.21	&	35.04	&	35.09	&	29.79	&	22.08	&	21.84	&	22.26	&	14.85 \\
\hline
t17	&	5.55	&	6.02	&	5.85	&	4.44	&	4.51	&	3.87	&	3.44	&	2.8	&	2.22	&	2.05 \\
\hline
t18	&	14.25	&	11.2	&	10.79	&	8.76	&	8.25	&	6.12	&	6.96	&	5.39	&	4.74	&	3.8 \\
\hline
t19	&	17.1	&	13.72	&	14.69	&	12	&	10.45	&	8.73	&	7.04	&	8.26	&	6.06	&	4.65 \\
\hline
t20	&	7.95	&	6.44	&	6.89	&	6.72	&	6.16	&	3.87	&	3.6	&	3.92	&	3.36	&	2.5 \\
\hline
t21	&	9.6	&	7.84	&	7.41	&	6.84	&	7.7	&	5.67	&	5.44	&	3.92	&	3.78	&	3.05 \\
\hline
t22	&	87.75	&	70.98	&	67.6	&	62.52	&	64.46	&	51.84	&	45.76	&	38.29	&	34.74	&	27.85 \\
\hline
t23	&	119.4	&	104.58	&	114.14	&	98.16	&	88.33	&	70.47	&	66.72	&	61.88	&	49.98	&	45.15 \\
\hline
t24	&	82.65	&	69.58	&	62.53	&	63.48	&	57.09	&	44.73	&	43.6	&	39.69	&	33.6	&	26.65 \\
\hline
\end{tabular}}
\label{WPP}
\end{center}
\end{table}

\begin{table}[h]
\begin{center}
\caption{Location of wind farms in the IEEE 3-Area RTS-96 system.}
 \scalebox{1}[1]{
\begin{tabular}{|c|c|c|c|c|c|c|c|c|c|}
\hline
 Unit &  Node & Unit &  Node &Unit &  Node &Unit &  Node &Unit &  Node \\
\hline\hline
$f_{1}$ & 103 & $f_{4}$ & 121 & $f_{7}$ & 216 & $f_{10}$ & 303 & $f_{13}$ & 316 \\
\hline
$f_{2}$ & 105 & $f_{5}$ & 203 & $f_{8}$ & 221  & $f_{11}$ & 305 & $f_{14}$ & 321 \\
\hline
$f_{3}$ & 116 & $f_{6}$ & 205 & $f_{9}$ & 223  & $f_{12}$ & 307 & $f_{15}$ & 323\\
\hline
\end{tabular}}
\label{wfloc}
\end{center}
\end{table}

\begin{table}[h]
\caption{Unit commitment plan (UCP), cost of commitment decisions (CCD) [\$], number of partitions, value of  weighting factor in Zhao model (ZWF), expected total cost (ETC) [\$] and worst-case total cost (WCTC) [\$] for the HUC and Zhao models applied to the IEEE 14-node system.}
 \begin{center}
 \scalebox{1}[1]{
\begin{tabular}{|c|c|c|c|c|c|}
\hline
UCP & CCD [\$] & \# Partitions & \# ZWF & ETC [\$] & WCTC [\$]
  \\
\hline
1 & 62569 & 6 - 10 & 0 - 0.2 & 286602 & 311534 \\
\hline
2 & 63505 & 1 - 5  & 0.2 - 1 & 287131 & 307030 \\
\hline
\end{tabular}}
\label{table:zhaopartitions} 
 \end{center}
\end{table}
%In order to prove the conservativeness of the commitment solution, Figure \ref{subfigure1} illustrates the results obtained in Figure \ref{CD Cost} when the expected total cost (ETC) and worst-case total cost (WCTC) performs respectively. We observe that higher number of partitions provide less conservative solutions and therefore, they adapt better when the ETC is applied and the opposite behavior occurs when a lower number of partitions is performed. 
%We observe that for a higher number of partitions, the WCTC increases while the ETC decreases and vice versa. The reason for this behavior can be found in how conservative is the commitment solution obtained in each partition. On the one hand, the solution committed for a lower number of partitions performs better in terms of cost when the worst-case total cost realizes, therefore they provide a more conservative solution. On the other hand, higher number of partitions provide less conservative solutions since they perform better when the expected cost is applied.

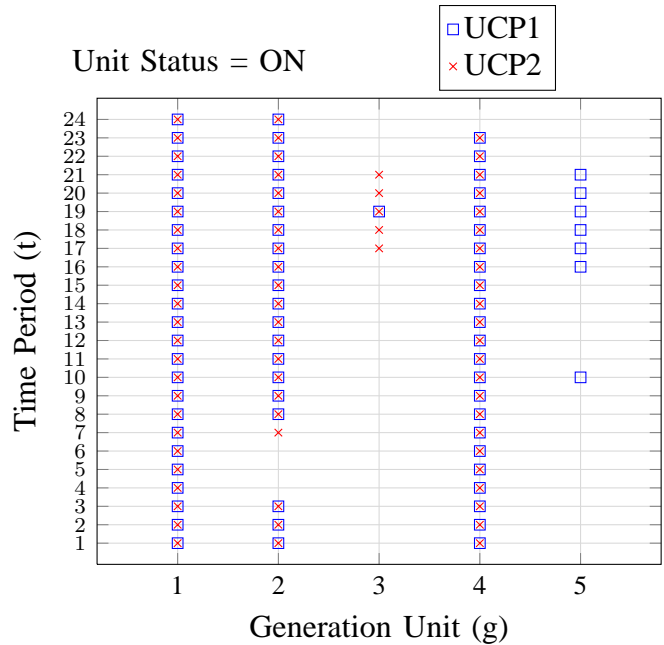
\begin{figure}[h]
\centering
\begin{tikzpicture}[scale=1]
\pgfplotsset{every x tick label/.append style={font=\small, yshift=0 ex}}
\pgfplotsset{every y tick label/.append style={font=\scriptsize, xshift=0 ex}}
	\begin{axis}[%
	grid=both,
	grid style={gray!30},
	 title style={xshift=-2.5cm,},
	legend style={at={(0.71,1.2)},anchor=north},
	title={Unit Status = ON},
	 xlabel=Generation Unit (g),
    ylabel=Time Period (t),
	symbolic x coords={1,2,3,4,5},
    xtick=data,
    ytick=data,
    enlarge x limits=0.2,
    enlarge y limits=0.05,
	scatter/classes={%
		a={mark=square,blue},%
		b={mark=x,red},%
		c={mark=o,draw=black}}]
	\addplot[scatter,only marks,%
		scatter src=explicit symbolic]%
	table[meta=label] {
x   y   label
1	1	a
1	2	a
1	3	a
1	4	a
1	5	a
1	6	a
1	7	a
1	8	a
1	9	a
1	10	a
1	11	a
1	12	a
1	13	a
1	14	a
1	15	a
1	16	a
1	17	a
1	18	a
1	19	a
1	20	a
1	21	a
1	22	a
1	23	a
1	24	a
2	1	a
2	2	a
2	3	a
2	8	a
2	9	a
2	10	a
2	11	a
2	12	a
2	13	a
2	14	a
2	15	a
2	16	a
2	17	a
2	18	a
2	19	a
2	20	a
2	21	a
2	22	a
2	23	a
2	24	a
3	19	a
4	1	a
4	2	a
4	3	a
4	4	a
4	5	a
4	6	a
4	7	a
4	8	a
4	9	a
4	10	a
4	11	a
4	12	a
4	13	a
4	14	a
4	15	a
4	16	a
4	17	a
4	18	a
4	19	a
4	20	a
4	21	a
4	22	a
4	23	a
5	10	a
5	16	a
5	17	a
5	18	a
5	19	a
5	20	a
5	21	a
};	\legend{}

	\addplot[scatter,only marks,%
		scatter src=explicit symbolic]%
	table[meta=label] {
x   y   label
1	1	b
1	2	b
1	3	b
1	4	b
1	5	b
1	6	b
1	7	b
1	8	b
1	9	b
1	10	b
1	11	b
1	12	b
1	13	b
1	14	b
1	15	b
1	16	b
1	17	b
1	18	b
1	19	b
1	20	b
1	21	b
1	22	b
1	23	b
1	24	b
2	1	b
2	2	b
2	3	b
2	7	b
2	8	b
2	9	b
2	10	b
2	11	b
2	12	b
2	13	b
2	14	b
2	15	b
2	16	b
2	17	b
2	18	b
2	19	b
2	20	b
2	21	b
2	22	b
2	23	b
2	24	b
3	17	b
3	18	b
3	19	b
3	20	b
3	21	b
4	1	b
4	2	b
4	3	b
4	4	b
4	5	b
4	6	b
4	7	b
4	8	b
4	9	b
4	10	b
4	11	b
4	12	b
4	13	b
4	14	b
4	15	b
4	16	b
4	17	b
4	18	b
4	19	b
4	20	b
4	21	b
4	22	b
4	23	b
};	\legend{UCP1 , UCP2}
	\end{axis}
\end{tikzpicture}
\caption{Unit-commitment plans obtained for the IEEE 14 Nodes system.}
\label{Figure14nodes}
\end{figure}

We can see from Table \ref{table:zhaopartitions} that UCP2 is more conservative than UCP1. Essentially, the main difference between both plans is that, as shown in Fig.~\ref{Figure14nodes}, UCP1 relies more on generating unit 5 than on unit 3. The former is comparatively smaller, but also more flexible. Therefore, the system can take advantage of unit 5 if the eventual wind power production turns out to be high.

%{\color{red} UCP1 relies more on generating unit 5 than on unit 3. The former is comparatively smaller and more expensive, but also more flexible. Therefore, the system can take advantage of unit 5 if the eventual wind power production turns out to be high.}

%the two commitment plans obtained behave differently for the expected cost and worst-case cost realization. These results are obtained fixing the commitment decision obtained for the previous cases in the SUC and RUC respectively. In brief, this proves that UCP1 is less conservative than UCP2 since this first one provides a better result in terms of expected cost than the the second one and vice versa. To see more in detail this assumption,  depicts both unit commitments. 

%If we take a look to the unit commitment plans, we can observe how UCP1 dispatch unit 5 while UCP2 dispatch unit 3 for time periods between 17 and 21. For those specific time periods the demand is the highest during the whole time horizon and the wind power production specially low. Therefore, both plans dispatch extra units during those time periods. Nevertheless, unit 3 (dispatch by UCP2), is a non flexible unit with more capacity than unit 5 (dispatch by UCP1) which is a flexible generator with less capacity. According to this planning production, UCP2 is adapting to a more pessimistic realization of the uncertainties while UCP1, by scheduling a more flexible unit is adapting to a wider realization of the uncertainties. As a consequence, the solution provided by UCP2 is more conservative than the solution provided for UCP1. Towards to asses this concept, 

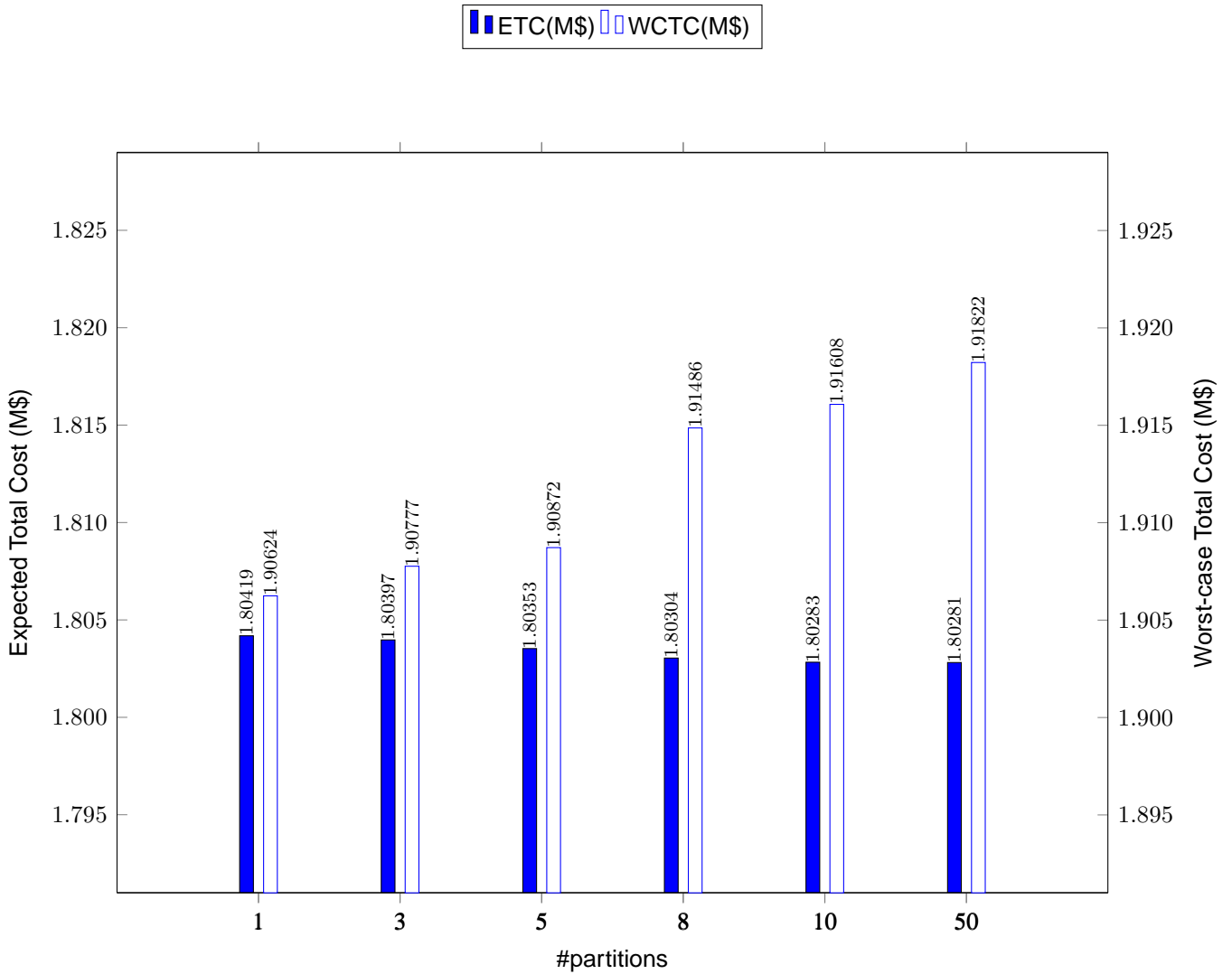
\begin{figure}[h]
 \centering
\begin{tikzpicture}[font=\sffamily\small,scale=1]
\pgfkeys{
    /pgf/number format/precision=5, 
    /pgf/number format/fixed zerofill=true
}

\begin{axis}[
  yticklabel style={/pgf/number format/precision=3},
  axis y line*=left,
  ymin=1.80, ymax=1.82,
  ybar, /pgf/bar shift=-5pt,
  enlargelimits=0.2,
  enlarge y limits=0.45,
  bar width = .2cm,
  legend style={at={(0.5,1.15)},
  anchor=north,legend columns=-1},
  xlabel={\#partitions},
  xticklabel style={/pgf/number format/precision=0,font=\small},
  ylabel={Expected Total Cost (M\$)},
  yticklabel style={/pgf/number format/precision=3},
  width=0.9\linewidth,
  height=0.7\linewidth,
  symbolic x coords={1,3,5,8,10,50},
  xtick=data,
  nodes near coords,
  nodes near coords align={vertical},
  every node near coord/.append style={color=black, rotate=90, anchor=center, font=\scriptsize, xshift=14, yshift=0}
]
\addplot [black,fill=blue]
  coordinates{
    (1,1.804186)
    (3,1.803973)
    (5,1.803529)
    (8,1.803043)
    (10,1.802829)
    (50,1.802813)
}; \label{plot_one} 
\end{axis}

\begin{axis}[
  yticklabel style={/pgf/number format/precision=3},
  axis y line*=right,
  ymin=1.90, ymax=1.92,
  ybar, /pgf/bar shift=5pt,
  enlargelimits=0.2,
  enlarge y limits=0.45,
  bar width = .2cm,
  legend style={at={(0.5,1.20)},
  anchor=north,legend columns=-1},
  xticklabel style={/pgf/number format/precision=0,font=\small},
  ylabel={Worst-case Total Cost (M\$)},
  yticklabel style={/pgf/number format/precision=3},
  width=0.9\linewidth,
  height=0.7\linewidth,
  symbolic x coords={1,3,5,8,10,50},
  xtick=data,
  nodes near coords,
  nodes near coords align={vertical},
  every node near coord/.append style={color=black, rotate=90, anchor=center, font=\scriptsize, xshift=14, yshift=0}
]
\addlegendimage{/pgfplots/refstyle=plot_one}\addlegendentry{ETC(M\$)}
\addplot[blue,fill=white]
  coordinates{
    (1,1.906244)
    (3,1.907766)
    (5,1.908722)
    (8,1.914864)
    (10,1.916075)
    (50,1.918224)
}; \label{plot_two} \addlegendentry{WCTC(M\$)}
\end{axis}
\end{tikzpicture}%}
      \caption{Expected total cost (ETC) and worst-case total cost (WCTC) of the IEEE 3-Area RTS-96 system for 1, 3, 5, 8, 10 and 50 partitions.}
        \label{subfigure1}
\end{figure}

Fig.~\ref{subfigure1} illustrates the expected and the worst-case system operating cost (denoted by ETC and WCTC, respectively) prompted by the HUC solution for a different number of partitions when the IEEE 3-Area RTS-96 system is solved. It is clear that, as we increase the number of partitions, the proposed HUC formulation provides less conservative unit-commitment plans, which gradually perform better in expectation but worse under the worst-case scenario. Indeed, the HUC solutions for 1 and 50 partitions correspond to the robust and the stochastic unit-commitment solutions, respectively. 

\subsection{Performance Evaluation of SPDA}\label{ResultsPerformance}
%

%%%%%%%%%%%%%%%%%%%%%%%%%%%%%%%%%%%%%%%%%%%%%%%%%%%%%%%%%%%%%%%%%%%%%%%%%%
%%%%%%                       NEW PART                     %%%%%%%%%%%%%%
%I should emphasize that with the algorithm I can achieve better results because I can reduce the MIP gap due to the simplification of the problem now it is tractable. Of course we have to say that it is a simplification. The algorithm simplifies the HUC and yields the same solution. Also the algorithm can perform better solutions, since the MIP gap can be reduced more bacause the problem is more simplified. How fast is our algortihm in percentage.

We now assess the effectiveness of the proposed parallelization-and-decomposition solution scheme on the IEEE 3-Area RTS-96 system. For this purpose, we solve the HUC model (\ref{46})--(\ref{76}) both directly and by means of SPDA. Thus, we compare these two alternative solution approaches, which are both coded in GAMS using CPLEX 12.6.1 and implemented in the DTU HPC Cluster. The DTU HPC Cluster is a composite of a variety of hardware components of different technical characteristics. Therefore, we refer the reader to \cite{Centr57:online} for further and detailed information on the cluster and its components. We solve the raw HUC (without decomposition) using one node in a multi-threading configuration that counts on 10 Intel cores, while SPDA is implemented in a multi-node and multi-threaded environment. In particular, partitions are solved in parallel, each in a different node of the cluster with up to 10 Intel cores per node. Lastly, the reduced HUC problem is solved using again one node employing up to 10 Intel cores.    

\begin{table}[h]
\caption{Solution time [seconds], number of variables, number of equations, number of scenarios finally retained, and time saving [\%]  for a different number of partitions.}
\label{BCWCTAB2}
\centering
 \scalebox{1}[1]{
\begin{tabular}{c c c c c c c}
\hline
\# Partitions &  & T [s] & \# VAR & \# EQU & \# SC &   $\Delta$T [\%]  \\ 
\hline
\multirow{2}{*}{1} &  HUC & 77172 & 287714 & 892273 & 50 &  \multirow{2}{*}{99.81} \\
\hhline{~-----~}
                     &  SPDA & 144 & 5977 & 17665 & 3 &  \\
\hline\hline
\multirow{2}{*}{3} &  HUC & 57184 & 287716 & 892273 & 50 &  \multirow{2}{*}{99.23}  \\
\hhline{~-----~}
                     &  SPDA & 439 & 46228 & 132678 & 7 & \\
\hline\hline
\multirow{2}{*}{5} &  HUC & 33591 & 287718 & 892273 & 50 &  \multirow{2}{*}{92.15}  \\
\hhline{~-----~}
                     &  SPDA & 2636 & 79926 & 238668 & 13  & \\
\hline\hline
\multirow{2}{*}{8} &  HUC &  41485 & 287721 & 892273 & 50 &  \multirow{2}{*}{92.30}  \\
\hhline{~-----~}
                     &  SPDA & 3192 & 113625 & 344658 &19 &  \\
\hline\hline
\multirow{2}{*}{10} &  HUC & 16458 & 287723 & 892273 & 50 & \multirow{2}{*}{82.68}  \\
\hhline{~-----~}
                     &  SPDA & 2849 & 136091 & 415318 & 23 & \\
\hline\hline
\multirow{2}{*}{50} &  HUC & \multirow{2}{*}{8346} & \multirow{2}{*}{287713} & \multirow{2}{*}{892223} & \multirow{2}{*}{50} &  \multirow{2}{*}{0}  \\
\hhline{~-~~~~~}
                     &  SPDA &  &  &  & &\\
\hline
\end{tabular}}
\end{table}

The performance of the aforementioned solution strategies are compared in Table \ref{BCWCTAB2}. The solutions yielded by SPDA are the same as or even better than the solutions provided by the raw HUC.  This is so, because SPDA ends up solving a smaller version of the raw HUC problem (\ref{46})--(\ref{76}) and, thereby, generally results in unit-commitment solutions with a lower MIP gap. Furthermore, SPDA drastically diminishes the required solution time when the number of partitions is kept low enough, in which case the number of scenarios considered in the last step of the algorithm is small  (see the last two columns of Table \ref{BCWCTAB2}). Needless to say, SPDA reduces to solving the raw HUC problem for a number of partitions equal to 50. Note that the number of binary variables to be treated is independent of the considered solution strategy and equal to 6912 in all cases.      

\section{Conclusion and Future Research}
\label{Conclusion}
In this paper we propose a new formulation of the unit commitment problem under uncertainty that allows us to find unit-commitment plans that perform relatively well in terms of both the expected and the worst-case system operating cost. The new formulation relies on clustering the scenario data set into a number of partitions. The expectation of the system operating cost is then taken over those scenarios that result in the worst-case dispatch cost within each partition. The conservatism of the so-obtained unit-commitment solution (that is, how close it is to the pure scenario-based stochastic or robust unit-commitment plan) is controlled via the user-specified number of partitions. 
We also develop a parallelization-and-decomposition scheme to efficiently solve the proposed unit-commitment formulation. Our numerical results show that our scheme is able to dramatically reduce the required running time while improving the optimality of the found solution. 

We envision two possible avenues of future research at least. First, we would like to explore the possibility of further using decomposition to solve the reduced version of the proposed unit commitment formulation that is obtained after applying our parallelization-and-decomposition scheme. Second, we would like to investigate how to extend our formulation and the associated solution approach to a multi-stage setup.

%, stochastic programming and robust optimization are merged in a hybrid stochastic-robust unit commitment for the integration of large amount of wind power production scenarios. Partitions over the complete set of scenarios are performed in order to merged stochastic and robust optimization. The number of partitions determines how conservative is the problem. Therefore, the purpose is to control the conservativeness of the final solution by varying the number of partitions. Furthermore, a scenario reduction technique is elaborated in order to improve the tractability and the computational efficiency of the hybrid stochastic-robust unit commitment problem. This scenario reduction technique, based on a Benders decomposition approach, leads to a double parallelization process which speeds up the solution time in a very efficient way. Future research may be focused on elaborate a method to define an optimal number of partitions according to the decision maker requirements.   

%\section*{Acknowledgment}


\ifCLASSOPTIONcaptionsoff
  \newpage
\fi

%Esto es todo


\begin{thebibliography}{1}

\bibitem{birge2011introduction}
J.~R. Birge and F.~Louveaux, {\em Introduction to stochastic programming}.
\newblock Springer Science \& Business Media, 2011.

\bibitem{ben2009robust}
A.~Ben-Tal, L.~El~Ghaoui, and A.~Nemirovski, {\em Robust optimization}.
\newblock Princeton University Press, 2009.


%\bibitem{barth2006stochastic}
%R.~Barth, H.~Brand, P.~Meibom, and C.~Weber, ``A stochastic unit-commitment
%  model for the evaluation of the impacts of integration of large amounts of
%  intermittent wind power,'' in {\em Probabilistic Methods Applied to Power
%  Systems, 2006. PMAPS 2006. International Conference on}, pp.~1--8, IEEE,
%  2006.


  \bibitem{takriti1996stochastic}
S.~Takriti, J.~R. Birge, and E. Long, ``A stochastic model for the unit commitment problem,'' {\em IEEE Transactions on Power Systems}, vol.~11, no.~3, pp.~1497--1508, 1996.

  \bibitem{bouffard2005market}
F.~Bouffard, F.~D. Galiana, and A.~J. Conejo, ``Market-clearing with stochastic
  security-part i: formulation,'' {\em IEEE Transactions on Power Systems},
  vol.~20, no.~4, pp.~1818--1826, 2005.

\bibitem{wang2008security} %PARA LOS SCENARIOS
J.~Wang, M.~Shahidehpour, and Z.~Li, ``Security-constrained unit commitment
  with volatile wind power generation,'' {\em IEEE Transactions on Power Systems}, vol.~23, no.~3, pp.~1319--1327, 2008.
  
\bibitem{morales2009economic}
J.~M. Morales, A.~J. Conejo, and J. P{\'e}rez-Ruiz, ``Economic valuation of reserves in power systems with high penetration of wind power,'' {\em IEEE Transactions on Power Systems},  vol.~24, no.~2, pp.~900--910, 2009.

\bibitem{tuohy2009unit}
A.~Tuohy, P.~Meibom, E.~Denny, and M.~O'Malley, ``Unit commitment for systems
  with significant wind penetration,'' {\em IEEE Transactions on Power Systems}, vol.~24, no.~2, pp.~592--601, 2009.
  
   \bibitem{wang2012chance}
Q.~Wang, Y.~Guan, and J.~Wang, ``A chance-constrained two-stage stochastic
  program for unit commitment with uncertain wind power output,'' {\em IEEE Transactions on Power Systems}, vol.~27, no.~1, pp.~206--215, 2012.
  
%    \bibitem{atamturk2007two}
%A.~Atamt{\"u}rk and M.~Zhang, ``Two-stage robust network flow and design under
%  demand uncertainty,'' {\em Operations Research}, vol.~55, no.~4,
%  pp.~662--673, 2007.
  
%  \bibitem{bertsimas2003robust}
%D.~Bertsimas and M.~Sim, ``Robust discrete optimization and network flows,''
%  {\em Mathematical programming}, vol.~98, no.~1, pp.~49--71, 2003.
  
    \bibitem{zeng2013solving}
B.~Zeng and L.~Zhao, ``Solving two-stage robust optimization problems using a
  column-and-constraint generation method,'' {\em Operations Research Letters},
  vol.~41, no.~5, pp.~457--461, 2013.
  
  \bibitem{zhao2012robust}
L.~Zhao and B.~Zeng, ``Robust unit commitment problem with demand response and
  wind energy,'' in {\em Power and Energy Society General Meeting, 2012 IEEE},
  pp.~1--8, IEEE, 2012. 
  
  \bibitem{bertsimas2013adaptive}
D.~Bertsimas, E.~Litvinov, X.~A. Sun, J.~Zhao, and T.~Zheng, ``Adaptive robust
  optimization for the security constrained unit commitment problem,'' {\em
  IEEE Transactions on Power Systems}, vol.~28, no.~1, pp.~52--63, 2013.
  
  \bibitem{Heitsch} H.~Heitsch and W.~R\"{o}misch, ``Scenario reduction algorithms in
stochastic programming,'' \emph{Computational Optimization and
Applications}, vol.~24, pp. 187-206, 2003.

\bibitem{Morales2} J.~M. Morales, S. Pineda, A.~J. Conejo, and M. Carri\'on,
``Scenario Reduction for Futures Market Trading in Electricity
Markets,'' \emph{IEEE Transactions on Power Systems}, vol.~24, no. 2, pp.
878--888, May 2009.

\bibitem{Salva} S. Pineda and A.~J. Conejo,
``Scenario reduction for risk-averse electricity trading,'' in \emph{IET Generation, Transmission \& Distribution}, vol.~4, no. 6, pp.
694--705, 2010.
  
%\bibitem{zugno2013robust}
%M.~Zugno and A.~J. Conejo, ``A robust optimization approach to energy and
%  reserve dispatch in electricity markets,'' tech. rep., Technical University
%  of Denmark, 2013.
  
  \bibitem{liu2003fuzzy}
L.~Liu, G.~Huang, Y.~Liu, G.~Fuller, and G.~Zeng, ``A fuzzy-stochastic robust
  programming model for regional air quality management under uncertainty,''
  {\em Engineering Optimization}, vol.~35, no.~2, pp.~177--199, 2003.

\bibitem{xu2009srccp}
Y.~Xu, G.~Huang, X.~Qin, and M.~Cao, ``SRCCP: a stochastic robust
  chance-constrained programming model for municipal solid waste management
  under uncertainty,'' {\em Resources, Conservation and Recycling}, vol.~53,
  no.~6, pp.~352--363, 2009.
  
    
  \bibitem{liubidding}
G.~Liu, Y.~Xu, and K.~Tomsovic, ``Bidding strategy for microgrid in day-ahead
  market based on hybrid stochastic/robust optimization,'' {\em
  IEEE Transactions on Smart Grid}, vol.~7, no.~1, pp.~227--237, 2016.
  
  
\bibitem{fanzeres2015contracting}
B.~Fanzeres, A.~Street, and L.~A. Barroso, ``Contracting strategies for
  renewable generators: a hybrid stochastic and robust optimization approach,''
  {\em IEEE Transactions on Power Systems}, vol.~30, no.~4, pp.~1825--1837,
  2015.
        
   \bibitem{zhao2013unified}
C.~Zhao and Y.~Guan, ``Unified stochastic and robust unit commitment,'' {\em
  IEEE Transactions on Power Systems}, vol.~28, no.~3, pp.~3353--3361, 2013.
  
\bibitem{zhao2015data}
C.~Zhao and Y.~Guan, ``Data-driven stochastic unit commitment for integrating
  wind generation,'' {\em IEEE Transactions on Power Systems}, vol.~31, no.~4, pp.~2587--2596, 2016.


\bibitem{Gabrel2014471}
V.~Gabrel and C.~Murat and A.~Thiele, ``Recent advances in robust optimization: {A}n overview,'' {\em European Journal of Operational Research}, vol.~235, no.~3, 2014.

\bibitem{wagstaff2001constrained}
K.~Wagstaff, C.~Cardie, S.~Rogers, S.~Schr{\"o}dl, {\em et~al.}, ``Constrained
  k-means clustering with background knowledge,'' in {\em ICML}, vol.~1,
  pp.~577--584, 2001.

\bibitem{hobbs1999stochastic}
B.~F. Hobbs and Y.~Ji, ``Stochastic programming-based bounding of expected
  production costs for multiarea electric power system,'' {\em Operations
  Research}, vol.~47, no.~6, pp.~836--848, 1999.

    \bibitem{watson2014new}
J.-P. Watson, F.~D. Munoz, and B.~Hobbs, ``New bounding and decomposition
  approaches for {MILP} investment problems: Multi-area transmission and
  generation planning under policy constraints.,'' tech. rep., Sandia National
  Laboratories (SNL-NM), Albuquerque, NM (United States), 2014.
  
    \bibitem{Papavasiliou}
A. Papavasiliou and S. S. Oren, ``Multi-Area Stochastic Unit Commitment for High Wind Penetration in a Transmission Constrained Network,'' {\em Operations Research}, vol.~61, no.~3, pp.~578--592, 2013.
  
%      \bibitem{morales2013integrating}
%J.~M. Morales, A.~J. Conejo, H.~Madsen, P.~Pinson, and M.~Zugno, {\em
%  Integrating renewables in electricity markets: operational problems},
%  vol.~205.
%\newblock Springer Science \& Business Media, 2013.  
  
  \bibitem{carrion2006computationally}
M.~Carri{\'o}n and J.~M. Arroyo, ``A computationally efficient mixed-integer
  linear formulation for the thermal unit commitment problem,'' {\em IEEE Transactions on Power Systems}, vol.~21, no.~3, pp.~1371--1378, 2006.
  
%    \bibitem{bussieck2009grid}
%M.~R. Bussieck, M.~C. Ferris, and A.~Meeraus, ``Grid-enabled optimization with
%  gams,'' {\em INFORMS Journal on Computing}, vol.~21, no.~3, pp.~349--362,
%  2009.

  \bibitem{Papavasiliou2015}
A.~Papavasiliou, S.~S.~Oren, B.~Rountree, ``Applying High Performance Computing to Transmission-Constrained Stochastic Unit Commitment for Renewable Penetration,'' {\em IEEE Transactions on Power Systems}, vol.~30, no.~3, pp.~1690--1701, 2015.

  \bibitem{Centr57:online}
``Central DTU HPC cluster |.'' \url{http://www.cc.dtu.dk/?page_id=342}.
\newblock (Accessed on 03/21/2016).

  \bibitem{GAMSD96:online}
``Gams documentation center.''
  \url{http://www.gams.com/help/index.jsp?topic=%2Fgams.doc%2Fsolvers%2Fguss%2Findex.html&anchor=GUSS_APP_DEA_MODELING}.
\newblock (Accessed on 06/01/2016).


\bibitem{Power25:online}
``Power Systems Test Case Archive - UWEE.''
  \url{https://www.ee.washington.edu/research/pstca/}.
\newblock (Accessed on 04/19/2016).

  \bibitem{TheRE28:online}
``The Real Lab - Renewable Energy Analysis Laboratory.''
  \url{http://www.ee.washington.edu/research/real/library.html}.
\newblock (Accessed on 03/21/2016).

    \bibitem{pinson2013wind}
P.~Pinson {\em et~al.}, ``Wind energy: Forecasting challenges for its
  operational management,'' {\em Statistical Science}, vol.~28, no.~4,
  pp.~564--585, 2013.
  






\end{thebibliography}
\end{document}